\documentclass[a4paper,11pt,reqno]{article}
\usepackage[dvipsnames]{xcolor}
\usepackage{amsmath}
\usepackage{amsthm}
\usepackage{amssymb}
\usepackage{mathcomp}
\usepackage{dsfont}
\usepackage{comment}
\usepackage{authblk}

\usepackage[active]{srcltx}
\usepackage{stackengine}
\usepackage[bookmarksnumbered=true]{hyperref}
\usepackage{bbm}

\usepackage{geometry}
\newtheorem{definition}{Definition}[section]

\newtheorem{proposition}[definition]{Proposition}
\newtheorem{theorem}[definition]{Theorem}
\newtheorem{remark}[definition]{Remark}

\numberwithin{equation}{section}

\usepackage[deletedmarkup=xout,authormarkup=none,
]{changes}

\begin{document}

\title{On the magnetic laplacian with a piecewise constant magnetic field in $\mathbb{R}^3_+$}
\author{Emanuela L. Giacomelli}
\affil{LMU M\"unich, Department of Mathematics, Theresienstr. 39, 80333 M\"unchen, Germany}
%
%
\maketitle

\abstract{We consider the Neumann realization of the magnetic laplacian in $\mathbb{R}^3_+$, in the case in which the magnetic field has a piecewise constant strength and a uniform direction. This operator is expected to be an effective model in studying the threshold to the normal state for a 3D superconductor exposed to a discontinuous magnetic field. We review some recent results above the infimum of the spectrum of the aforementioned operator.}

\section{Introduction}\label{sec:int}

We study a Schr\"odinger operator  in $\mathbb{R}^3_+:= \{\mathbf{x}\in \mathbb{R}^3\, \vert\, \mathbf{x} = (x_1, x_2, x_3), \, x_3>0\}$ with a magnetic field admitting a piecewise constant strength and a uniform direction. We review some recent results obtained in \cite{AG} about the bottom of the spectrum. Such an operator is interesting to be considered in the theory of superconductivity. We first introduce it and later we motivate the last assertion.

\subsection{The main result} We denote by $\mathbf{B}_{\alpha,\gamma,a}$ a piecewise constant magnetic field having uniform direction. More precisely, we suppose $\mathbf{B}_{\alpha,\gamma,\alpha}$ to have an intensity equal to $1$ in the region $\mathcal{D}^1_\alpha$ and equal to $a$ in $\mathcal{D}^2_\alpha$, where
\[
    \mathcal{D}^1_\alpha = \big\{\mathbb{R}^3\ni\mathbf{x} = \rho(\cos\theta\sin\phi, \sin \theta\sin \phi, \cos\phi)\, \vert \, \rho\in (0,\infty), \theta\in (0,\alpha), \, \phi\in (0,\pi)\big\},
\]
\[
   \mathcal{D}^1_\alpha = \big\{\mathbb{R}^3\ni\mathbf{x} = \rho(\cos\theta\sin\phi, \sin \theta\sin \phi, \cos\phi)\, \vert \, \rho\in (0,\infty), \theta\in (\alpha,\pi), \, \phi\in (0,\pi)\big\}.
\]
For any  $a\in [-1, 1)\setminus \{0\}$, $\alpha\in (0,\pi)$ and $\gamma \in [0,\pi/2]$, we define 
\begin{equation}
    \mathbf{B}_{\alpha, \gamma,a} = (\cos\alpha \sin\gamma, \sin\alpha\sin\gamma, \cos\gamma)\big(\mathbbm{1}_{\mathcal{D}^1_\alpha} + a\mathbbm{1}_{\mathcal{D}^2_\alpha}\big).
\end{equation}
Note that by symmetry considerations, we can restrict the angle $\gamma$ to the interval $[0,\pi/2]$. Moreover, instead of taking into account a generic piecewise constant magnetic field having two different uniform intensities in the regions $\mathcal{D}^1_\alpha$ and $\mathcal{D}^2_\alpha$, it is possible (by scaling arguments) to reduce the study to the case of a field having intensity equal to $\mathbbm{1}_{\mathcal{D}^1_\alpha} + a \mathbbm{1}_{\mathcal{D}^2_\alpha}$ for some $a\in [-1,1)$. The condition $a\neq0$ we imposed above is a technical restriction which we need to simplify our analysis.

Let now $\mathbf{A}_{\alpha,\gamma,a}\in H^{1}_{\mathrm{loc}}(\mathbb{R}^3_+, \mathbb{R}^3)$ be a vector potential such that $\mathrm{curl}\mathbf{A}_{\alpha,\gamma,a} = \mathbf{B}_{\alpha,\gamma,a}$. In this paper we take into account the magnetic Neumann realization of the following self-adjoint operator in $\mathbb{R}^3_+$
\begin{equation}\label{eq: def op L}
    \mathcal{L}_{\alpha, \gamma,a} = -(\nabla - i \mathbf{A}_{\alpha,\gamma,a})^2,
\end{equation}
with a domain given by 
\begin{multline}
    \mathrm{Dom}(\mathcal{L}_{\alpha,\gamma,a}) = \big\{ u\in L^2(\mathbb{R}^3_+)\quad : \quad (\nabla - i \mathbf{A}_{\alpha,\gamma,a})^n u\in L^2(\mathbb{R}^3_+), 
    \\
    \mbox{for}\,\,\, n\in \{1,2\}, \,\, (\nabla - i \mathbf{A}_{\alpha,\gamma,a})u \cdot (0,1,0)\vert_{\partial{\mathbb{R}^3_+}} = 0\big\}.
\end{multline} 
We denote by $\lambda_{\alpha,\gamma,a}$ the bottom of the spectrum of $\mathcal{L}_{\alpha,\gamma,a}$, i.e., 
\begin{equation}\label{eq: def lambda}
    \lambda_{\alpha,\gamma,a} := \inf \mathrm{sp}(\mathcal{L}_{\alpha,\gamma,a}).
\end{equation}

The main result proved in \cite{AG} is in the theorem below.
\begin{theorem}[Bottom of the spectrum of $\mathcal{L}_{\alpha,\gamma,a}$]\label{thm: main} Let $a\in [-1,1)\setminus \{0\}$, $\alpha\in (0,\pi)$, $\gamma \in [0,\pi/2]$ and $\nu_0= \mathrm{arcsin}(\sin\alpha\sin\gamma)$. Let $\lambda_{\alpha,\gamma,a}$ as in \eqref{eq: def lambda}. It holds
\begin{equation}\label{eq: ineq lambda main}
    \lambda_{\alpha,\gamma,a} \leq \mathrm{min}(\beta_a, |a|\zeta_{\nu_0}),
\end{equation}
where $\beta_a$ and $\zeta_{\nu_0}$ are as in \eqref{eq: beta a mu a} and \eqref{eq: def zeta nu}.
Moreover, if the inequality in \eqref{eq: ineq lambda main} is strict, then there exists $\tau_\star\in\mathbb{R}$ such that 
\[
    \lambda_{\alpha,\gamma,a} = \underline{\sigma}(\alpha,\gamma, a, \tau_\star),
\]
where $\underline{\sigma}(\alpha,\gamma,a, \tau_\star)$ is an eigenvalue of the operator $\underline{\mathcal{L}}_{\alpha,\gamma,a}(\tau_\star)$ defined in \eqref{eq: def underline L}. 
\end{theorem}

\subsection{Motivation}

The motivations to study the operator $\mathcal{L}_{\alpha,\gamma,a}$ go back to the phenomenon of superconductivity, which was discovered in 1911 by H. Kamerlingh Onnes in Leiden. Here we put the focus on the breakdown of superconductivity in presence of an external magnetic field \cite{saint1963onset}. In general, superconductors can be divided into two types, according to how the breakdown occurs. For type-I, superconductivity is abruptly destroyed via a first order phase transition. In 1957 Abrikosov deduced the existence of a class of materials which
exhibit a different behavior, i.e., some of their superconducting properties are preserved
when submitted to a suitably large magnetic field. Physically, these two classes can be
identified by the value of a parameter $\kappa$, also known as the Ginzburg-Landau parameter (i.e., a value proportional to the inverse of the penetration depth and typical of the material). The value $\kappa$ is smaller than $1/\sqrt{2}$ for type-I superconductors and larger than $1/\sqrt{2}$ for the so called type-II superconductors.
We consider here extreme type-II superconductors, i.e., we assume that the
Ginzburg-Landau parameter satisfies the condition $\kappa\gg 1$. 

In general, it is well-known \cite{giorgi2002breakdown} that a superconducting material exposed to a strong magnetic field with intensity $h_{\mathrm{ex}}$ loses permanently its superconducting properties (i.e., goes to the normale state) when $h_{\mathrm{ex}}$ exceeds some critical value.  Determining this critical value is not an easy task and strongly depends on the geometry of the sample. Below we underline the main ideas toward this characterization, and emphasizing where the operator we introduced in \eqref{eq: def op L} is expected to play a role.
{}
\subsubsection{The Ginzburg-Landau theory}
To study the transition to the normal state, it is convenient to use the Ginzburg-Landau theory \cite{GL} which, in general, allows to describe the behavior of a type-II superconductor exposed to an external magnetic field $\mathbf{B}$ (such that $|\mathbf{B}| = h_{\mathrm{ex}}$) in a temperature close to the critical one. Let $\Omega \subset\mathbb{R}^3$ be a domain, the GL functional is defined by\footnote{Note that in some cases the GL model reduces to a two dimensional one: this happens when for example the external magnetic field is supposed to be perpendicular to the cross section of a superconducting wire.}
\[
\mathcal{G}_{\Omega, \kappa}[\psi,\mathbf{A}] = \int_\Omega\, |\nabla + i \mathrm{h}_{\mathrm{ex}}\mathbf{A})\psi|^2 - \frac{1}{2}\kappa^2(2|\psi|^2 - |\psi|^4) + h_{\mathrm{ex}}\int_{\mathbb{R}^3} |\mathrm{curl}\mathbf{A} - 1|^2.
\]
Here $\psi$ is the order parameter ($|\psi|^2$ denotes the density of Cooper pairs) and $h_{\mathrm{ex}}\mathbf{A}$ is the induced magnetic vector potential. We recall that we take into account extreme type II-superconductors, which means that $\kappa \rightarrow \infty$. 

Using the GL theory it is possible to study the phase transitions which occur in a type-II superconductor: these can be described  by identifying three increasing critical values of the magnetic field. When the first critical value $H_{c_1}$ is reached, superconductivity is lost in the bulk of the sample at isolated points (see e.g.\cite{SS,fournais2010spectral}). Between the second and third critical fields, i.e., in the regime $H_{c_2} \leq h_{\mathrm{ex}} \leq H_{c_3}$, superconductivity survives only close to the boundary of the sample (see e.g.,\cite{Pan1, CR1,CR2, CG1, HK, CG2,CG3}) . Above the third critical field $H_{c_3}$, the sample goes back to its normal state (see e.g., \cite{lu99eigen, helffer2001magnetic, HP, fournais2006third, FH07, BNF}).

In this framework the normal state corresponds to the choice $(\psi, \mathbf{A}) = (0, \mathbf{F})$ where $\mathbf{F}$ is such that $\mathrm{curl}\mathbf{F} = 1$, i.e., three are no Cooper pairs and the external field penetrates completely the sample. It is then natural to expect that to characterize the value of $H_{c_3}$, the first term in GL functional is playing the main role, i.e.,  one should study the magnetic Laplacian:
\[
    -(\nabla - i h_{\mathrm{ex}}\mathbf{A})^2.
\]
As suggested in \cite[Chapter 13]{fournais2010spectral}, we take into account external magnetic fields of intensity proportional to the GL parameter $\kappa$, i.e., $h_{\mathrm{ex}} = \kappa \sigma $ for some $\sigma>0$. Note that since $\kappa\rightarrow \infty$, the intensity of the magnetic field is high. Under this choices, the study of $H_{c_3}$ is naturally linked to a semiclassical limit $h\rightarrow 0$:
\begin{equation}\label{eq: semiclassical}
   -(\nabla - i(\kappa\sigma)\mathbf{A})^2 = - (\kappa\sigma)^2(h\nabla - i \mathbf{A})^2, \qquad h:= (\kappa\sigma)^{-1}.
\end{equation}
\subsubsection{A semiclassical problem}
Analysing the  semiclassical problem in \eqref{eq: semiclassical} is important to prove the localization of the GL minimizing order parameter in order to characterize the transition to the normal state. Many works indeed have been dedicated to the study of the operator $(-ih\nabla - \mathbf{A})^2$ in the limit $h\rightarrow 0$ deriving an asymptotics of the first eigenvalue and proving the localization of the associated eigenfunctions. In dimensions $d=2,3$, the first eigenvalue of $(h\nabla - i \mathbf{A})^2$ behaves at first order as $h\mathcal{E}(\mathbf{B}, \Omega)$, where $\mathcal{E}(\mathbf{B}, \Omega)$ is the smallest eigenvalue of a given model operator which strongly depends on the dimension $d$, on the geometry of $\Omega$ and on the shape of the magnetic field. In other words, to study the semiclassical problem it is useful to first take into account specific effective models.

Below we list some well-know situations and we underline where the 3D Schr\"odinger operator we study here is expected to appear.
In the case of a uniform external magnetic, if the domain $\Omega\subset\mathbb{R}^d$ ($d=2,3$) is smooth, one has to deal with two model operators: one defined over $\mathbb{R}^d$ (when working in the interior of $\Omega$) and the other living in $\mathbb{R}^d_+$ (to work near the boundary $\partial\Omega$). For such situations we refer to \cite{lu99eigen,helffer2001magnetic,fournais2010spectral} for dimension $d=2$ and to  \cite{giorgi2002breakdown,lu2000surface,helffer2004magnetic,fournais2006third, fournais2010spectral} in the 3D setting. Moreover, it is well-know that in presence of singularities (e.g, corners, wedges) along the boundary of $\Omega$, one has to introduce another operator to work close to the singularities. This operator is defined over an infinite angular sector or an infinite wedges (according to the dimension), see \cite{ BNF,bonnaillie2015ground} for the 2D case and \cite{ lu2000surface, Pan2,
 popoff2013schrodinger,bonnaillie2015ground, popoff2015model} for $d=3$. For not-uniform magnetic field there are, in general, less results available in the literature. We mention here \cite{Assaad3, assaad2020magnetic,AssaadHearing21} for piecewise constant magnetic fields in dimension $d=2$ and \cite{raymond2009sharp, raymond2013var} for smoothly varying magnetic fields both in dimension $d=2$ and $d=3$. 

The transition to the normal state in the case of a smooth domain $\Omega\subset\mathbb{R}^3$ with an external magnetic field which is piecewise constant is completely open. The study of the 3D Schr\"odinger operator introduced above can be seen as a first step towards this characterization. To give an idea of that, we take into account an external magnetic such that 
\[
    \mathbf{B}(x) = (\mathbbm{1}_{x_2>0} + a\mathbbm{1}_{x_2<0})(x)(0,0,1)\qquad \mbox{for}\,\, x= (x_1, x_2,x_3),
\]
and a domain $\Omega\subset\mathbb{R}^3$ with smooth boundary which intersects transversally the plane $(x_1x_3)$, we call this intersection the \textit{discontinuity plane} and we denote it by $S$. Moreover, we set $\Gamma:= \partial\Omega \cap S$ to be the discontinuity curve. In this setting superconductivity is expected to nucleate close to $\Gamma$ right before disappearing. To  rigorously prove this, it is convenient to use the model operator $\mathcal{L}_{\alpha,\gamma,a}$ introduced above to work in regions localized at the boundary along $\Gamma$. Indeed, the fact that localization occurs at the boundary force us to work with a model operator defined on $\mathbb{R}^3_+$. Moreover, to work close to $\Gamma$ requires to take into account a model operator with a discontinuous magnetic field. \\

\noindent\textbf{Acknowledgments.} The author acknowledges the support of the Istituto Nazionale di Alta Matematica ``F. Severi", through the Intensive Period ``INdAM Quantum Meetings (IQM22)".

\section{Proof of Theorem \ref{thm: main}}
In this section we give the main ideas for the study of the bottom of the spectrum of 
\[
    \mathcal{L}_{\alpha,\gamma,a} = -(\nabla - i \mathbf{A}_{\alpha,\gamma,a}), \qquad \alpha \in (0,\pi), \gamma\in\left[0,\frac{\pi}{2}\right],\quad a\in [-1,1)\setminus\{0\},
\]
i.e., $\lambda_{\alpha, \gamma,a} = \inf\mathrm{sp}(\mathcal{L}_{\alpha,\gamma,a})$. In the following we make a specific choice of the vector potential $\mathbf{A}_{\alpha,\gamma,a}$ which allows us to use a partial Fourier transform to compare our operator with 2D operators.

\subsection{The reference 2D operators}

As mentioned above, it turns out that it is possible to compare $\mathcal{L}_{\alpha,\gamma,a}$ with 2D operators. We will study the spectrum of such 2D operators making use of two additional models. The first one is a 2D operator with discontinuous magnetic field and the second one is Schr\"odinger operator in $\mathbb{R}^3_+$ having a uniform magnetic field. We introduce them in what follows and we refer to \cite{AG} (and references therein) for more details.

\subsubsection{Magnetic laplacian with a piecewise constant magnetic field in 2D}\label{subsec: beta a mu a} Let $\mathbf{A}_a\in H^1_{\mathrm{loc}}(\mathbb{R}^2, \mathbb{R}^2)$ be such that 
\[
    \mathrm{curl}\mathbf{A}_a(x) = \big(\mathbbm{1}_{x_2 >0} + a\mathbbm{1}_{x_2 <0}\big)(x), \qquad x\in \mathbb{R}^2, \quad a\in [-1,1)\setminus\{0\}.
\]
Consider the magnetic Neumann realization of 
\begin{equation}
    \mathcal{L}_a:= \big(-\nabla - i \mathbf{A}_a\big)^2,
\end{equation}
with domain 
\[
    \mathrm{dom}(\mathcal{L}_a) := \big\{u\in L^2(\mathbb{R}^2)\quad :\quad (\nabla - i\mathbf{A}_a)^nu\in L^2(\mathbb{R}^2), \,\,n=1,2\big\}.
\]
We denote by $\beta_a$ the bottom of the spectrum, i.e., 
\begin{equation}
    \beta_a =\inf\mathrm{sp}(\mathcal{L}_a).
\end{equation}
The operator $\mathcal{L}_a$ as well as the value $\beta_a$ were widely studied (see e.g., \cite{Assaadlowest20, Assaad2019, hislop2016band, hislop2015edge}). Here we just recall that $\mathcal{L}_a$ can decomposed by one dimensional fiber operator via a partial Fourier transform, i.e.,
\begin{equation}
    \mathcal{L}_a = \int_{\xi\in\mathbb{R}}^{\oplus}\mathfrak{h}_a(\xi) \, d\xi, \qquad \mathfrak{h}_a(\xi) = \begin{cases} -\frac{d^2}{dt^2} + (t-\xi)^2,  &\mbox{for}\,\,\, t>0 \\ 
    -\frac{d^2}{dt^2} + a(t-\xi)^2,  &\mbox{for}\,\,\, t<0.
    \end{cases}
\end{equation}
As a consequence, denoting by $\mu_a(\xi)$ the bottom of the spectrum of $\mathfrak{h}_a(\xi)$, one has 
\begin{equation}\label{eq: beta a mu a}
    \beta_a = \inf_{\xi\in\mathbb{R}}\mu_a(\xi).
\end{equation}
See \cite[Section 2.1]{AG} for the main properties of $\beta_a$ and $\mu_a(\cdot)$.

\subsubsection{Magnetic laplacian with constant magnetic field in $\mathbb{R}^3_+$} 

We now take into account a uniform magnetic field $\mathbf{B}_\nu$ with unit strength on $\mathbb{R}^3_+$, where $\nu$ denotes the angle between $\mathbf{B}_\nu$ and the plane $(x_1x_3)$. We can explicitly write 
\[
    \mathbf{B}_\nu = (0,\sin\gamma,\cos\gamma).
\]
We can then take into account the magnetic Neumann realization of 
\begin{equation}
    H_\nu = -(\nabla - i\mathbf{A}_\nu)^2 \qquad\mbox{in}\,\,\, L^2(\mathbb{R}^3_+),
\end{equation}
where $\mathbf{A}_\nu\in H^1_{\mathrm{loc}}(\mathbb{R}^3_+, \mathbb{R}^3)$ is such that $\mathrm{curl}\mathbf{A}_\nu = \mathbf{B}_\nu$. We denote the bottom of the spectrum of $H_\nu$ by 
\begin{equation}\label{eq: def zeta nu}
    \zeta_\nu := \inf\mathrm{sp}(H_\nu).
\end{equation}
This model operator  is studied in  \cite{lu99eigen, lu2000surface, morame2005remarks}. We refer to \cite{AG} for a collection of some useful properties of $\zeta_\nu$.

\subsection{Ideas for the proof of the main result}

We now summarize the strategy of the proof of Theorem \ref{thm: main} done in \cite{AG}. More precisely, first we reduce the study of $\lambda_{\alpha,\gamma,a}$ to the one of the bottom of the spectrum of 2D operators (Section \ref{sec: reduction}), then we collect the main properties we need on the spectrum of the aforementioned 2D operators (Section \ref{sec: spectrum 2D}) and in Section \ref{sec: conclusion} we give an idea of the final proof. 

\subsubsection{Reduction to 2D operators}\label{sec: reduction}
Now we want to do a partial Fourier transform to decompose (see Section \ref{sec: reduction}) $\mathcal{L}_{\alpha,\gamma,a}$. To do that, it is convenient to fix the gauge. Thus, from now on we suppose that the vector potential $\mathbf{A}_{\alpha,\gamma,a}$ is such that $\mathbf{A}_{\alpha,\gamma,a} = (A_1, A_2, A_3)$, with 
\begin{eqnarray*}
    A_1 &=& 0, \\
    A_2&=& \begin{cases}  \cos\gamma\, \big(x_1 - x_2\,(1-a) \cot\alpha\big) &\mbox{for}\,\,\, x\in \mathcal{D}_\alpha^1,\\ a\cos\gamma\, x_1  &\mbox{for}\,\,\, x\in \mathcal{D}_\alpha^2, \end{cases} \\
    A_3&=& \begin{cases} \sin\gamma \, \big( x_2\cos\alpha - x_1\sin\alpha \big) &\mbox{for}\,\,\, x\in \mathcal{D}_\alpha^1, \\
    a\sin\gamma \, \big( x_2\cos\alpha - x_1\sin\alpha \big) &\mbox{for}\,\,\, x\in \mathcal{D}_\alpha^2.
     \end{cases}
\end{eqnarray*}
Note that this choices for $A_1$, $A_2$, $A_3$ ensure that $\mathbf{A}_{\alpha,\gamma,a} \in H^{1}_{\mathrm{loc}}(\mathbb{R}^3_+, \mathbb{R}^3)$ and imply that the operator $\mathcal{L}_{\alpha,\gamma,a}$ is translation invariant with respect to the $x_3$ coordinate. We can then use a partial Fourier transform in the $x_3$ variable to decompose $\mathcal{L}_{\alpha,\gamma,a}$ via fiber operators living in $\mathbb{R}_+^2$. More precisely, we can write 
\begin{equation}\label{eq: decomposition fiber}
    \mathcal{L}_{\alpha,\gamma,a} = \int_{\tau\in\mathbb{R}}^{\oplus}\, \underline{\mathcal{L}}_{\alpha,\gamma,a}(\tau)\, d\tau, 
\end{equation}
where 
\begin{equation}\label{eq: def underline L}
    \underline{\mathcal{L}}_{\alpha,\gamma,a}(\tau) = -(\nabla - i \underline{\mathbf{A}}_{\alpha,\gamma,a})^2 + V_{\alpha,\gamma,a}(\tau)
\end{equation}
Below we explain our notations. First, we set $D^1_\alpha$, $D^2_\alpha$ to be the orthogonal projections of the regions $\mathcal{D}^1_\alpha$, $\mathcal{D}^2_\alpha$ over the plane $(x_1x_2)$. The magnetic potential $\underline{\mathbf{A}}_{\alpha,\gamma,a}$ is the projection of $\mathbf{A}_{\alpha,\gamma,a}$ on $\mathbb{R}^2_+$, i.e., $\underline{\mathbf{A}}_{\alpha,\gamma,a} = (\underline{A}_1, \underline{A}_2)$ with $\underline{A}_1=0$ and
\begin{equation}
 \underline{A}_2 = \begin{cases} \cos\gamma (x_1 - (1-a)\cot\alpha x_2) &\mbox{for}\,\,\, (x_1, x_2)\in D^1_\alpha  \\ a\cos\gamma &\mbox{for}\,\,\, (x_1,x_2)\in D^2_\alpha.\end{cases} 
\end{equation}
Moreover, $\underline{\mathbf{A}}_{\alpha,\gamma,a}$ is such that 
\begin{equation}
    \mathrm{curl}\underline{\mathbf{A}}_{\alpha,\gamma,a} = \underline{s}_{\alpha,a}\cos\gamma, \qquad \underline{s}_{\alpha,a} = \mathbbm{1}_{D^1_\alpha} + a \mathbbm{1}_{D^2_\alpha}.
\end{equation}
Finally the potential $V_{\alpha,\gamma,a}(\tau)$ appearing in \eqref{eq: def underline L} is an electric potential which is defined through the projection of $\mathbf{B}_{\alpha,\gamma,a}$ on $\mathbb{R}^2_+$. More precisely, we denote the aforementioned projection by $\underline{\mathbf{B}}_{\alpha,\gamma,a}$ and, explicitly, we have 
\begin{equation}
    \underline{\mathbf{B}}_{\alpha,\gamma,a} = (\cos\alpha, \sin\alpha\sin\gamma)\underline{\mathbf{s}}_{\alpha,a}\equiv (\underline{b}_1, \underline{b}_2).
\end{equation}
The electric potential is then given by
\begin{equation}
    V_{\alpha,\gamma,a}(\tau) = (x_1 \underline{b}_2 - x_2 \underline{b}_1 -\tau)^2.
\end{equation}
From \eqref{eq: decomposition fiber} it turns out that 
\[
    \lambda_{\alpha,\gamma,a} = \inf_\tau\underline{\sigma}_{\alpha,\gamma,a}(\tau),
\]
where we denoted by $\underline{\sigma}_{\alpha,\gamma,a}(\tau)$ the bottom of the spectrum of the operator $\underline{\mathcal{L}}_{\alpha,\gamma,a}(\tau)$.
As a consequence, we reduced the study of $\lambda_{\alpha,\gamma,a}$ to the one of the map $\tau\mapsto\underline{\sigma}(\alpha,\gamma,a)$ (which can be proven to be $C^\infty$).

\subsubsection{Spectrum of the 2D operators}\label{sec: spectrum 2D}

Here we recall two results we need about the spectrum of the 2D reduced operator $\underline{\mathcal{L}}_{\alpha,\gamma,a}(\tau)$ for fixed\footnote{The case $\gamma = 0$ can be treated directly, this is why suppose $\gamma \neq 0$ in this section.} $\alpha\in(0,\pi)$, $\gamma\in(0,\pi/2]$, $a\in [-1,1)\setminus\{0\}$, $\tau\in\mathbb{R}$. 

\begin{proposition}[Bottom of the essential spectrum] \label{pro: essential spectrum 2D} Let $a \in [-1,1)\setminus\{0\}$, $\alpha\in (0,\pi)$, $\gamma\in (0, \pi/2]$ and $\tau\in\mathbb{R}$. Let
\begin{equation}    
    \underline{\sigma}_{\mathrm{ess}}(\alpha,\gamma,a,\tau) = \inf\mathrm{sp}_{ess}(\underline{\mathcal{L}}_{\alpha,\gamma,a}).
\end{equation}
It holds 
\begin{equation}
    \underline{\sigma}_{ess}(\alpha,\gamma,a,\tau) = \inf_{\xi\in\mathbb{R}}\big(\mu_a(\tau\sin\gamma + \xi\cos\gamma) + (\xi\sin\gamma - \tau\cos\gamma)^2\big),
\end{equation}
where $\mu_a(\cdot)$ is as in \eqref{eq: beta a mu a}.
\end{proposition}
\begin{proposition}[Behavior of $\underline{\sigma}(\alpha,\gamma,a,\tau)$ for large $\tau$] \label{pro: spectrum limit tau }Let $\alpha\in (0,\pi)$ and $\gamma\in (0,\pi/2]$. It holds:
\begin{enumerate}
\item For $a\in [-1,0)$:
\[
    \lim_{\tau\rightarrow - \infty} \underline{\sigma}(\alpha,\gamma,a,\tau) = +\infty, \qquad \lim_{\tau\rightarrow +\infty}\underline{\sigma}(\alpha,\gamma,a,\tau) = |a| \zeta_{\nu_0}.
\]
\item For $a\in (0,1)$, 
\[
    \lim_{\tau\rightarrow - \infty} \underline{\sigma}(\alpha,\gamma,a,\tau) = a\zeta_{\nu_0}, \qquad \lim_{\tau\rightarrow +\infty}\underline{\sigma}(\alpha,\gamma,a,\tau) = \zeta_{\nu_0}.
\]
\end{enumerate}
\end{proposition}
The proofs of Proposition \ref{pro: essential spectrum 2D} and Proposition \ref{pro: spectrum limit tau } are based on the study of the two auxiliary operators: one is useful to work near the boundary of $\mathbb{R}^2_+$ away from the discontinuity line (i.e., the intersection between the plane of equation $x_1\sin\alpha - x_2\cos\alpha = 0$ and $\mathbb{R}^2_+$), meanwhile the other is an effective operator useful when working close to the discontinuity. We refer to \cite[Section 3]{AG} for more details.

\subsubsection{Conclusion of the proof of Theorem \ref{thm: main}}\label{sec: conclusion} 

Once Proposition \ref{pro: essential spectrum 2D} and Proposition \ref{pro: spectrum limit tau } are established, the proof of Theorem \ref{thm: main} is quite simple. As mentioned before, one can distinguish between $\gamma = 0$ and $\gamma \neq 0$. In the first case, it is immediate to get that 
\begin{equation}\label{eq: gamma 0 1}
    \lambda_{\alpha,0,a} \leq |a| \Theta_0,
\end{equation}
by following what was proven in \cite[Section 3]{Assaad3}. Combining \eqref{eq: gamma 0 1} with the fact that $\zeta_0 = \Theta_0$ and that $\beta_a \geq |a| \Theta_0$, one has
\begin{equation}
    \lambda_{\alpha, 0, a} \leq \min (\beta_a, |a| \zeta_0).
\end{equation}
We refer to \cite[Section 4]{AG} for more details. 

In the case $\gamma\neq 0$ the proof is more involved. In particular, from Proposition \ref{pro: spectrum limit tau }, we get that 
\begin{equation}\label{eq: est 1 final}
    \underline{\sigma}(\alpha, \gamma,a, \tau) \leq |a| \zeta_{\nu_0}.
\end{equation}
We can now distinguish between $a\in [-1,0)$ and $a\in (0,1)$. In the second case, i.e., $a\in (0,1)$, there is nothing to prove. Indeed, one has that $\beta_a= a$ for $a\in(0,1)$ and that $\zeta_{\nu_0} < 1$ (see \cite[Section 2]{AG}). This allows to conclude the proof of \eqref{eq: ineq lambda main} for $a$ positive. 
In the case $\gamma\neq 0$, $a\in [-1,0)$, we have to work a bit more. From Proposition \ref{pro: essential spectrum 2D} and choosing a particular value\footnote{One has to take $\tau_\ast = \xi_a\sin\gamma$, where $\xi_a$ is the minimum of $\mu_a(\cdot)$ introduced in \eqref{eq: beta a mu a}.} of $\tau = \tau_\ast$, one has 
\begin{equation}\label{eq: est 2 final}
    \underline{\sigma}(\alpha,\gamma,a,\tau_\ast) \leq \beta_a.
\end{equation}
Combining \eqref{eq: est 2 final} with \eqref{eq: est 1 final}, the estimate in \eqref{eq: ineq lambda main} holds.
We now discuss the case of a strict inequality. From Proposition \ref{pro: spectrum limit tau }, we have 
\begin{equation}
\inf_{\tau} \underline{\sigma}(\alpha,\gamma,a, \tau) = \lambda_{\alpha,\gamma,a} < |a| \zeta_{\nu_0} = \min \big(\lim_{\tau\rightarrow -\infty} \underline{\sigma}(\alpha,\gamma,a,\tau), \lim_{\tau\rightarrow+\infty}\underline{\sigma}(\alpha,\gamma,a,\tau)\big),
\end{equation}
which implies that $\inf_{\tau} \underline{\sigma}(\alpha,\gamma,a, \tau)$ is attained at some $\tau_\star\in\mathbb{R}$. Moreover, from Proposition \ref{pro: essential spectrum 2D} (see \cite[Corollary 3.6]{AG}) we know that 
\[
    \inf_{\tau\in\mathbb{R}}\underline{\sigma}_{ess}(\alpha,\gamma,a,\tau) \geq \beta_a.
\]
Thus, we get 
\begin{equation}
    \lambda_{\alpha,\gamma,a} = \underline{\sigma}(\alpha,\gamma,a,\tau_\star) < \beta_a \leq \underline{\sigma}_{ess} (\alpha,\gamma,a,\tau_\star),
\end{equation}
which implies that $\lambda_{\alpha,\gamma,a}$ is an eigenvalue of $\underline{\mathcal{L}}_{\alpha,\gamma,a}(\tau_\star)$.

\begin{remark}\label{rem: 1}
In \cite[Proposition 1.4]{AG} we also provide a condition on $(\alpha,\gamma,a)$ such that the strict inequality in \eqref{eq: ineq lambda main} is realized.
\end{remark}

\begin{remark}
We consider cases of $(\alpha,\gamma,a,\tau)$ where the infimum of the spectrum of $\underline{\mathcal L}_{\alpha,\gamma,a}(\tau)$ is an eigenvalue below  the essential spectrum (see Remark \ref{rem: 1}).
One can prove an Agmon-estimate result showing the decay of the corresponding eigenfunction, for large values of $|x|$. More precisely, let $a\in[-1,1)\setminus\{0\}$,  $\alpha\in(0,\pi)$, $\gamma\in(0,\pi/2]$ and $\tau\in \mathbb{R}$. Consider the case where $\underline\sigma(\alpha,\gamma,a,\tau)<\underline\sigma_{ess}(\alpha,\gamma,a,\tau)$. Let $v_{\alpha,\gamma,a,\tau}$ be the  normalized eigenfunction corresponding to $\underline\sigma(\alpha,\gamma,a,\tau)$. For all $\eta\in \sqrt{\underline\sigma_{ess}(\alpha,\gamma,a,\tau)-\underline\sigma(\alpha,\gamma,a,\tau)}$, there exists a constant $C$  such that
    \[\underline Q_{\alpha,\gamma,a}^\tau(e^{\eta\phi}v_{\alpha,\gamma,a,\tau})\leq C,\]
    where  $\phi(x)=|x|$, for $x\in\mathbb{R}^2_+$ and $\underline Q_{\alpha,\gamma,a}^\tau$ is the quadratic form associated to $\underline{\mathcal{L}}_{\alpha,\gamma,a}(\tau)$ in $\mathbb{R}^2_+$. For the proof,  we refer the reader to similar results in~\cite[Theorem~9.1]{bonnaillie2003analyse} and \cite{bonnaillie2012discrete}.
\end{remark}

\end{document}